\documentclass[psamsfonts,reqno]{amsart}
\usepackage{amssymb,eucal}

\theoremstyle{plain}

\newtheorem{lemma}{Lemma}[section]
\newtheorem{theorem}[lemma]{Theorem}
\newtheorem{proposition}[lemma]{Proposition}
\newtheorem{corollary}[lemma]{Corollary}

\newtheorem*{stat}{\name}
\newcommand{\name}{testing}

\theoremstyle{definition}
\newtheorem{definition}[lemma]{Definition}

\theoremstyle{remark}

\newtheorem{notation}[lemma]{Notation}

\newenvironment{all}[1]{\renewcommand{\name}{#1}\begin{stat}}
                        {\end{stat}}

\newcommand{\qedc}{{\qed}~{\rm Claim~{\theclaim}.}}

\numberwithin{equation}{section}

\newcommand{\pup}[1]{\textup{(}{#1}\textup{)}}
\newcommand{\jirr}{join-ir\-re\-duc\-i\-ble}

\newcommand{\set}[1]{\{#1\}}
\newcommand{\setm}[2]{\set{#1\mid#2}}
\newcommand{\seq}[1]{\left\langle{#1}\right\rangle}
\newcommand{\seqm}[2]{\seq{{#1}\mid{#2}}}
\newcommand{\famm}[2]{\left(#1\mid#2\right)}
\newcommand{\bt}{\mathbin{\bowtie}}

\newcommand{\Pow}{\mathfrak{P}}
\newcommand{\dnw}{\mathbin{\downarrow}}

\DeclareMathOperator{\supp}{supp}
\DeclareMathOperator{\otp}{otp}
\DeclareMathOperator{\rk}{rk}
\newcommand{\cx}[1]{\Vert{#1}\Vert}
\newcommand{\cC}{\mathcal{C}}

\newcommand{\cR}{\mathcal{R}}

\newcommand{\cS}{\mathcal{S}}
\newcommand{\cD}{\mathcal{D}}
\newcommand{\cF}{\mathcal{F}}
\newcommand{\es}{\varnothing}

\newcommand{\onto}{\twoheadrightarrow}

\newcommand{\ba}{\boldsymbol{a}}
\newcommand{\bb}{\boldsymbol{b}}
\newcommand{\bc}{\boldsymbol{c}}
\newcommand{\bx}{\boldsymbol{x}}
\newcommand{\by}{\boldsymbol{y}}
\newcommand{\bz}{\boldsymbol{z}}
\newcommand{\bu}{\boldsymbol{u}}
\newcommand{\bv}{\boldsymbol{v}}
\newcommand{\bw}{\boldsymbol{w}}

\DeclareMathOperator{\Conc}{Con_c}

\newcommand{\ol}[1]{\overline{{#1}}}

\newcommand{\jz}{$\langle\vee,0\rangle$}

\newcommand{\jzs}{\jz-semi\-lat\-tice}

\newcommand{\jzh}{\jz-ho\-mo\-mor\-phism}

\newcommand{\jze}{\jz-em\-bed\-ding}

\newcommand{\jh}{join-ho\-mo\-mor\-phism}

\begin{document}

\title[Semilattice-valued measures on posets]%
{Non-extendability of semilattice-valued measures on partially ordered
sets}

\author[F.~Wehrung]{Friedrich Wehrung}
\address{CNRS, UMR 6139\\
D\'epartement de Math\'ematiques, BP 5186\\
Universit\'e de Caen, Campus 2\\
14032 Caen cedex\\
France}
\email{wehrung@math.unicaen.fr}
\urladdr{http://www.math.unicaen.fr/\~{}wehrung}

\keywords{Semilattice, poset, distributive, isotone, measure,
$\Delta$-Lemma, closed unbounded}
\subjclass[2000]{06A12, 06A06, 06A05}

\date{\today}

\begin{abstract}
For a poset $P$ and a distributive \jzs\ $S$, a \emph{$S$-valued
poset measure} on~$P$ is a map $\mu\colon P\times P\to S$ such that
$\mu(x,z)\leq\mu(x,y)\vee\mu(y,z)$, and $x\leq y$ implies that
$\mu(x,y)=0$, for all $x,y,z\in P$. In relation with congruence lattice
representation problems, we consider the problem whether such a measure
can be extended to a poset measure $\ol{\mu}\colon\ol{P}\times\ol{P}\to
S$, for a larger poset~$\ol{P}$, such that for all $\ba,\bb\in S$ and all
$x\leq y$ in $\ol{P}$, $\ol{\mu}(y,x)=\ba\vee\bb$ implies that there are a
positive integer~$n$ and a decomposition $x=z_0\leq z_1\leq\cdots\leq
z_n=y$ in~$\ol{P}$ such that either $\ol{\mu}(z_{i+1},z_i)\leq\ba$ or
$\ol{\mu}(z_{i+1},z_i)\leq\bb$, for all $i<n$.

In this note we prove that this is not possible as a rule, even in case
the poset~$P$ we start with is a \emph{chain} and $S$ has size $\aleph_1$.
The proof uses a ``monotone refinement property''
that holds in $S$ provided $S$ is either a lattice, or countable, or
strongly distributive, but fails for our counterexample. This strongly
contrasts with the analogue problem for \emph{distances} on (discrete)
sets, which is known to have a positive (and even \emph{functorial})
solution.
\end{abstract}

\maketitle

\section{Introduction}\label{S:Intro}

In the paper \cite{FPLat}, the author proved that for any lattice $K$,
any distributive \emph{lattice}~$S$ with zero, and any \jzh\ $\varphi$
from the \jzs\ $\Conc K$ of all finitely generated congruences of $K$ to
$S$, there are a lattice $L$, a lattice homomorphism
$f\colon K\to L$, and an isomorphism $\alpha\colon\Conc L\to S$ such that
$\varphi=\alpha\circ\Conc f$. In the paper
\cite{TuWe1}, J. T\r{u}ma and the author proved that for a \jzs~$S$,
this statement characterizes $S$ being a lattice. The proof of this
negative result strongly uses the lattice structure of the hypothetical
lattice~$L$, see the proof of \cite[Corollary~1.3]{TuWe1}.

In the present paper, we show that for a certain semilattice $S$ of
cardinality~$\aleph_1$, the poset structure alone is sufficient to get a
related counterexample. More precisely, for a \jzs\ $S$, a
\emph{$S$-valued poset measure} on a poset $P$ is a map
$\mu\colon P\times P\to S$ such that $\mu(x,z)\leq\mu(x,y)\vee\mu(y,z)$
(\emph{triangular inequality}) and $x\leq y$ implies that
$\mu(x,y)=\nobreak 0$, for all $x,y,z\in P$. We say that $\mu$ is a
\emph{V-measure}, if for all $x\leq y$ in $P$ and all $\ba,\bb\in S$, if
$\mu(y,x)\leq\ba\vee\bb$, then there are a positive integer~$n$ and a
decomposition $x=z_0\leq z_1\leq\cdots\leq z_n=y$ in~$P$ such that
either $\mu(z_{i+1},z_i)\leq\ba$ or $\mu(z_{i+1},z_i)\leq\bb$ for all
$i<n$. In particular, if $P$ is a lattice and $S=\Conc P$, then the map $\mu$ defined by $\mu(x,y)=\Theta^+(x,y)=\Theta(y,x\vee y)$ is a $\Conc P$\ -valued V-measure on $P$.

This yields the following poset analogue of the abovementioned
lattice-theoretical problem.

\begin{all}{Problem}
Let $S$ be a distributive \jzs. Does any $S$-valued poset measure on a
given poset extend to some $S$-valued poset V-measure on a larger poset?
\end{all}

A version of this problem for so-called \emph{distances} (instead of
measures) on discrete sets (instead of posets) is stated in \cite{RTW}.
The answer to this related question turns out to be positive (and easy).
More surprisingly, this positive solution can be made \emph{functorial}.

Nevertheless, we prove in the present paper that the problem above has a
\emph{negative} solution. Unlike what is done in
\cite{TuWe1}, we do not reach here a characterization of all lattices
among distributive \jzs s. Our counterexample, denoted by
$\cF(\omega_1)$ (see Corollary~\ref{C:NoIncrSeq2}) is obtained as an
application of a certain ``free construction'' used by M.
Plo\v{s}\v{c}ica and J. T\r{u}ma in \cite{PlTu}. The semilattice~$D$ of
\cite[Section~2]{TuWe1}, which is the simplest example of a \jzs\ which
is not a lattice, does not satisfy the negative property used here. This
is because~$D$ is \emph{countable}, while we prove in
Proposition~\ref{P:MonRef} that no countable distributive \jzs\ can have
the required negative property. On the other hand, in relation to
\cite[Problem~4]{TuWe1}, the proof of our counterexample uses very little
of the Axiom of Choice (namely, only the Axiom of \emph{countable}
choices), while the proof of the negative property of the abovementioned
semilattice $D$ established in
\cite[Corollary~2.4]{TuWe1} uses the existence of an embedding from
$\omega_1$ into the reals.

\section{Basic concepts}\label{S:Basic}

For posets (i.e., partially ordered sets)
$P$ and $Q$, a map $f\colon P\to Q$ is \emph{isotone}, if
$x\leq y$ implies that $f(x)\leq f(y)$, for all $x,y\in P$. In
addition, we say that $f$ is \emph{join-preserving}, if for any subset
$X$ of $P$, whenever the join $\bigvee X$ of $X$ exists in $P$, $\bigvee
f[X]$ exists in~$Q$, and $\bigvee f[X]=f\bigl(\bigvee X\bigr)$. For a
subset $X$ of a poset $P$, we shall put
$\dnw X=\setm{p\in P}{\exists x\in X\text{ such that }p\leq x}$,
and then $\dnw a=\dnw\set{a}$, for all $a\in P$. We say that $X$
is a \emph{lower subset of $P$}, if $X=\dnw X$.

A \jzs\ $S$ is \emph{distributive}, if $\bc\leq\ba\vee\bb$ in $S$ implies
that there are $\bx\leq\ba$ and $\by\leq\bb$ in $S$ such that $\bc=\bx\vee\by$. A distributive \jzs\ $S$ is \emph{strongly distributive}, if every element of $S$ is the join of a finite set of \jirr\ elements of $S$; equivalently, $S$ is
isomorphic to the semilattice of all finitely generated lower subsets of
some poset.

We shall denote by $\otp{P}$ the order-type of a well-ordered set $P$.
Hence $\otp{P}$ is an ordinal. We shall also use standard set-theoretical
notation and terminology, referring the reader to \cite{Jech} for further
information. In particular, we shall denote by~$\omega_1$ the first
uncountable ordinal. A subset $C$ of $\omega_1$ is \emph{closed
unbounded}, if $C$ is unbounded in $\omega_1$ and the join of any
nonempty bounded subset of $C$ belongs to~$C$. It is well-known that the
closed unbounded subsets form a countably complete filterbasis
on~$\omega_1$, see \cite[Lemma~7.4]{Jech}. Hence containing a closed
unbounded set is a notion of ``largeness'' for subsets of~$\omega_1$.

\section{Free distributive extension of a \jzs}\label{S:FDE}

There are several non-equivalent definitions of what should be the ``free
distributive extension'' of a given \jzs. The one that we shall use is
introduced in \cite[Section~2]{PlTu}. Let us first recall the
construction.

For a \jzs\ $S$, we shall put
$\cC(S)=\setm{\seq{\bu,\bv,\bw}\in S^3}{\bw\leq\bu\vee\bv}$.
A \emph{finite} subset $\bx$ of $\cC(S)$ is \emph{reduced}, if it
satisfies the following conditions:
\begin{itemize}
\item[(1)] $\bx$ contains exactly one diagonal triple, that is, a triple
of the form $\seq{\bu,\bu,\bu}$; we put $\bu=\pi(\bx)$.

\item[(2)] $\seq{\bu,\bv,\bw}\in\bx$ and $\seq{\bv,\bu,\bw}\in\bx$ implies
that $\bu=\bv=\bw$, for all $\bu,\bv,\bw\in\nobreak S$.

\item[(3)]
$\seq{\bu,\bv,\bw}\in\bx\setminus\set{\seq{\pi(\bx),\pi(\bx),\pi(\bx)}}$
implies that $\bu,\bv,\bw\nleq\pi(\bx)$, for all $\bu,\bv,\bw\in S$.
\end{itemize}

We denote by $\cR(S)$ the set of all reduced subsets of $\cC(S)$, endowed
with the partial ordering $\leq$ defined by
 \begin{equation}\label{Eq:DefleqRed}
 \bx\leq\by\ \Longleftrightarrow\
\forall\seq{\bu,\bv,\bw}\in\bx\setminus\by,
 \text{ either }\bu\leq\pi(\by)\text{ or }\bw\leq\pi(\by).
 \end{equation}
Furthermore, we shall identify $\bx$ with the
element $\set{\seq{\bx,\bx,\bx}}$ of $\cR(S)$, for all $\bx\in S$. For
set-theoretical purists, this can for example be done by replacing
$\cR(S)$ by the disjoint union of $S$ with the set of non-singletons in
$\cR(S)$. The disjointness can easily be achieved by a suitable
modification of the standard definition of a triple. We shall use the
symbol $\bt$ to denote the canonical generators of $\cR(S)$, so that
 \begin{equation*}
 \bt(\bu,\bv,\bw)=\begin{cases}
 \bw,&\text{if either }\bu=\bv\text{ or }\bv=0\text{ or }\bw=0,\\
 0,&\text{if }\bu=0,\\
 \set{\seq{0,0,0},\seq{\bu,\bv,\bw}},&\text{otherwise},
 \end{cases}
 \end{equation*}Observe that the canonical map $\pi\colon\cR(S)\onto S$ is
\emph{isotone} and that the restriction of $\pi$ to $S$ is the
identity. Furthermore, $\pi(\bt(\bu,\bv,\bw))=0$, for any non-diagonal
$\seq{\bu,\bv,\bw}\in\cC(S)$. The following is an easy consequence of
\eqref{Eq:DefleqRed}.
 \begin{equation}\label{Eq:SrelcpleteR(S)}
 \bx\leq\by\ \Longleftrightarrow\ \bx\leq\pi(\by),
 \qquad\text{for all }\seq{\bx,\by}\in S\times\cR(S).
 \end{equation}
We recall the standard facts established in \cite{PlTu} about this
construction.

\begin{proposition}\label{P:StR(S)}\hfill
\begin{enumerate}
\item For any \jzs\ $S$, $\cR(S)$ is a \jzs, and the inclusion map from
$S$ into $\cR(S)$ is a \jze.

\item For \jzs s $S$ and $T$, every \jzh\ $f\colon S\to T$ extends to a
unique \jzh\ $\cR(f)\colon\cR(S)\to\cR(T)$ such that
$\cR(f)(\bt(\bu,\bv,\bw))=\bt(f(\bu),f(\bv),f(\bw))$, for all
$\seq{\bu,\bv,\bw}\in\cC(S)$.

\item The assignment $S\mapsto\cR(S)$, $f\mapsto\cR(f)$ is a functor.
\end{enumerate}
\end{proposition}

The extension $\cR(S)$ is defined in such a way that $\bt(\bu,\bv,\bw)\leq
\bu$ and $\bw=\bt(\bu,\bv,\bw)\vee\bt(\bv,\bu,\bw)$, for all
$\seq{\bu,\bv,\bw}\in\cC(S)$. Hence, putting $\cR^0(S)=S$ and
$\cR^{n+1}(S)=\cR(\cR^n(S))$ for each~$n$, we obtain that the increasing
union $\cD(S)=\bigcup\famm{\cR^n(S)}{n<\omega}$ is a distributive \jzs,
extending~$S$. Furthermore,
putting $\cD(f)=\bigcup\famm{\cD^n(f)}{n<\omega}$ for each \jzh\ $f$, we
obtain that~$\cD$ is a functor.
The proof of the following lemma is straightforward.

\begin{lemma}\label{L:RD(inters)}
Let $S$ be a \jzs\ and let $\seqm{S_i}{i\in I}$ be a family of
\jz-subsemilattices of $S$. The following statements hold:
\begin{enumerate}
\item $\cR\left(\bigcap_{i\in I}S_i\right)=\bigcap_{i\in I}\cR(S_i)$
and $\cD\left(\bigcap_{i\in I}S_i\right)=\bigcap_{i\in I}\cD(S_i)$.
\item If $I$ is a nonempty upward directed poset and $\seqm{S_i}{i\in I}$ is
isotone, then
$\cR\left(\bigcup_{i\in I}S_i\right)=\bigcup_{i\in I}\cR(S_i)$ and
$\cD\left(\bigcup_{i\in I}S_i\right)=\bigcup_{i\in I}\cD(S_i)$.
\end{enumerate}
\end{lemma}

\begin{definition}\label{D:Cplxx}
For a \jzs\ $S$ and an element $\bx\in\cD(S)$, we define the \emph{rank}
of $\bx$, denoted by $\rk\bx$, as the least natural number $n$ such that
$\bx\in\cR^n(S)$, and the \emph{complexity} of $\bx$, denoted by
$\cx{\bx}$, by $\cx{\bx}=0$ if $\bx\in S$, and
 \[
 \cx{\bx}=\sum\famm{\cx{\bu}+\cx{\bv}+\cx{\bw}+1}
 {\seq{\bu,\bv,\bw}\in\bx},
 \qquad\text{for all }\bx\in\cD(S)\setminus S.
 \]
\end{definition}

\section{The semilattices $\cS(\Lambda)$ and $\cF(\Lambda)$}\label{S:CSX}

For any \emph{chain} $\Lambda$, we shall denote by $\cS(\Lambda)$ the
\jzs\ defined by generators $\ba$, $\bb$, and $\bc_i$, for $i\in\Lambda$,
and relations $\bc_i\leq\ba\vee\bb$ and $\bc_i\leq\bc_j$, for all
$i\leq j$ in~$\Lambda$. Hence the elements of $\cS(\Lambda)$ either
belong to $\cS(\es)=\set{0,\ba,\bb,\ba\vee\bb}$ or have the form $\bc_i$,
$\ba\vee\bc_i$, or $\bb\vee\bc_i$, for some $i\in\Lambda$. We shall
identify $\cS(X)$ with the \jz-subsemilattice of $\cS(\Lambda)$ generated
by $\cS(\es)\cup\setm{\bc_i}{i\in X}$, for any
$X\subseteq\Lambda$.

For chains $X$ and $Y$, any isotone map $f\colon X\to Y$ gives
raise to a unique \jzh\ $\cS(f)\colon\cS(X)\to\cS(Y)$ fixing $\ba$ and
$\bb$ and sending $\bc_i$ to $\bc_{f(i)}$, for all $i\in X$. Of course,
the assignment $\Lambda\mapsto\cS(\Lambda)$, $f\mapsto\cS(f)$ is a
functor.

We denote by $\cF=\cD\circ\cS$ the composition of the two functors $\cD$
and $\cS$.

The proof of the following lemma is straightforward.

\begin{lemma}\label{L:S(inters)}
Let $\Lambda$ be a chain and let $\seqm{X_i}{i\in I}$ be a family of
subsets of $\Lambda$. The following statements hold:
\begin{enumerate}
\item $\cS\left(\bigcap_{i\in I}X_i\right)=\bigcap_{i\in I}\cS(X_i)$.
\item If $I$ is a nonempty upward directed poset and $\seqm{X_i}{i\in I}$ is
isotone, then
$\cS\left(\bigcup_{i\in I}X_i\right)=\bigcup_{i\in I}\cS(X_i)$.
\end{enumerate}
\end{lemma}

As an easy consequence of Lemmas~\ref{L:RD(inters)} and~\ref{L:S(inters)},
we get the following.

\begin{lemma}\label{L:Supp}
Let $\Lambda$ be a chain. Then for any $\bx\in\cF(\Lambda)$, there exists
a least \pup{with respect to the inclusion} subset $X$ of $\Lambda$ such
that $\bx\in\cF(X)$; this subset is finite. 
\end{lemma}

We denote by $\supp(\bx)$ the subset given by Lemma~\ref{L:Supp}, and we
call it the \emph{support} of~$\bx$.

\begin{notation}\label{Not:x[X/Y]}
For a chain $\Lambda$, well-ordered subsets $X$ and $Y$ of
$\Lambda$ such that $\otp{X}\leq\otp{Y}$, and
$\bx\in\cF(X)$, we set $\bx[Y/X]=\cF(e_{X,Y})(\bx)$, where $e_{X,Y}$
denotes the unique embedding from $X$ into $Y$ whose range is a lower
subset of $Y$.
\end{notation}

Hence $\bx[Y/X]$ belongs to $\cF(Y)$, for all $\bx\in\cF(X)$.

\begin{lemma}\label{L:Fix[Y/X]}
Let $\Lambda$ be a chain and let $X$, $Y$ be well-ordered subsets of
$\Lambda$ such that $\otp{X}\leq\otp{Y}$ and $X\cap Y$ is a lower subset
of both $X$ and $Y$. Then $\bx[Y/X]=\bx$, for all $\bx\in\cF(X\cap Y)$.
\end{lemma}

\begin{proof}
As the set $Z=X\cap Y$ is a lower subset of both $X$ and $Y$, the
homomorphism $\cF(e_{Z,X})$ (resp., $\cF(e_{Z,Y})$) is the inclusion map
from $\cF(Z)$ into $\cF(X)$ (resp., $\cF(Y)$). In particular,
$\bx=\cF(e_{Z,X})(\bx)=\cF(e_{Z,Y})(\bx)$. Therefore,
 \begin{equation*}
 \bx[Y/X]=\cF(e_{X,Y})(\bx)=\cF(e_{X,Y})\circ\cF(e_{Z,X})(\bx)
 =\cF(e_{Z,Y})(\bx)=\bx.\tag*{\qed}
 \end{equation*}
\renewcommand{\qed}{}
\end{proof}

We are now reaching a crucial lemma.

\begin{lemma}[Interpolation Lemma]\label{L:Interp}
Let $\Lambda$ be a chain, let $X$, $Y$ be finite subsets of~$\Lambda$,
and let $\seq{\bx,\by}\in\cF(X)\times\cF(Y)$. If $\bx\leq\by$, then either
there exists $\bz\in\cF(X\cap Y)$ such that $\bx\leq\bz\leq\by$ or
\pup{$Y\not\subseteq X$ and $\bc_{\min(Y\setminus X)}\leq\by$}.
\end{lemma}

\begin{proof}
We shall denote by $\pi^l_k$ the canonical map from $\cR^l\cS(\Lambda)$
onto $\cR^k\cS(\Lambda)$, for all natural numbers $k\leq l$.
Put $m=\rk\bx$ and $n=\rk\by$. Observe that $\supp(\bx)\subseteq X$ and
$\supp(\by)\subseteq Y$. We argue by induction on
$\cx{\bx}+\cx{\by}$. If either $\supp(\bx)\subseteq Y$ or
$\supp(\by)\subseteq X$ then either $\bz=\bx$ or $\bz=\by$ belongs to
$\cF(X\cap Y)$ and satisfies the inequalities $\bx\leq\bz\leq\by$, so we
are done. So suppose that $\supp(\bx)\not\subseteq Y$ and
$\supp(\by)\not\subseteq X$. In particular, $X\not\subseteq Y$ and
$Y\not\subseteq X$, and both $\supp(\bx)$ and $\supp(\by)$ are nonempty.
We put $\xi=\min(Y\setminus X)$.

Suppose that $m=n=0$, that is, $\bx,\by\in\cS(\Lambda)$. Pick
$i\in\supp(\bx)$. As\linebreak
$\bc_i\leq\bx\leq\by$, we obtain that either
$\by=\ba\vee\bb$ (a contradiction, as then $\supp(\by)=\es$) or
$\by\in\set{\bc_j,\ba\vee\bc_j,\bb\vee\bc_j}$ for some
$j\geq i$. If $i=j$, then $\supp(\bx)=\supp(\by)=\set{i}$, a
contradiction. If $i<j$, then $\xi=j$ and so $\bc_{\xi}\leq\by$.

Suppose now that $m<n$. Then $\bx\leq\by$ means that
$\bx\leq\pi^n_m(\by)$ (use \eqref{Eq:SrelcpleteR(S)}).
As $\pi^n_m(\by)$ has support contained in $Y$ and
rank at most $m$, it follows from the induction hypothesis that either
$\bc_{\xi}\leq\pi^n_m(\by)$ (thus, \emph{a fortiori}, $\bc_{\xi}\leq\by$)
or there exists $\bz\in\cF(X\cap Y)$ such that
$\bx\leq\bz\leq\pi^n_m(\by)$ (thus, \emph{a fortiori}, $\bz\leq\by$).

So suppose from now on that $m>0$ (i.e., $\bx\notin\cS(\Lambda)$) and
$m\geq n$. If $\bx=\bigvee_{i<k}\bx_i$ where $k\geq2$ and each~$\bx_i$
has support contained in $X$ and complexity less than $\cx{\bx}$, then we
apply the induction hypothesis to each inequality $\bx_i\leq\by$, for
$i<k$. If $\bc_{\xi}\nleq\by$, then for all $i<k$, there exists
$\bz_i\in\cF(X\cap Y)$ such that $\bx_i\leq\bz_i\leq\by$. Hence
$\bx\leq\bz\leq\by$, where $\bz=\bigvee_{i<k}\bz_i$ belongs to
$\cF(X\cap Y)$. This reduces the problem to the case where
$\bx=\bt(\bu,\bv,\bw)$, where $\seq{\bu,\bv,\bw}$ is a non-diagonal
triple of elements of $\cR^{m-1}\cS(X)$ of complexity less than
$\cx{\bx}$.

If $m>n$, then, as $\bt(\bu,\bv,\bw)=\bx\leq\by$ with
$\supp(\bx)\not\subseteq Y$,
$\bu,\bv,\bw\in\cR^{m-1}\cS(X)$, and $\by\in\cR^{m-1}\cS(Y)$, it
follows from \eqref{Eq:DefleqRed} that either $\bu\leq\by$ or
$\bw\leq\by$. If, for example, $\bu\leq\by$, then, by the induction
hypothesis, either $\bc_{\xi}\leq\by$ (in which case we are done) or there
exists $\bz\in\cF(X\cap Y)$ such that $\bu\leq\bz\leq\by$. In the second
case, $\bx\leq\bz\leq\by$. The argument is similar in case $\bw\leq\by$.

The remaining case is $m=n>0$. As $\bt(\bu,\bv,\bw)=\bx\leq\by$ with
$\supp(\bx)\not\subseteq\nobreak Y$,
$\bu,\bv,\bw\in\cR^{m-1}\cS(X)$, and $\by\in\cR^m\cS(Y)$,
it follows from \eqref{Eq:DefleqRed} that either
$\bu\leq\nobreak\pi^m_{m-1}(\by)$ or $\bw\leq\pi^m_{m-1}(\by)$. If
$\bu\leq\pi^m_{m-1}(\by)$, then, by the induction hypothesis, either
$\bc_{\xi}\leq\pi^m_{m-1}(\by)$ (thus, \emph{a fortiori},
$\bc_{\xi}\leq\by$) or there exists $\bz\in\cF(X\cap Y)$ such that
$\bu\leq\bz\leq\pi^m_{m-1}(\by)$ (in which case $\bx\leq\bz\leq\by$). The
case where $\bw\leq\pi^m_{m-1}(\by)$ is similar.
\end{proof}

\begin{lemma}\label{L:Supci}
Let $\Lambda$ be a chain and let $X$ be a nonempty subset of $\Lambda$
admitting a supremum, say, $\xi$, in $\Lambda$. Then $\bc_{\xi}$ is the
supremum of $\setm{\bc_i}{i\in X}$ in $\cF(\Lambda)$.
\end{lemma}

\begin{proof}
Let $\bx\in\cF(\Lambda)$ such that
$\bc_i\leq\bx$ for all $i\in X$, we prove that $\bc_{\xi}\leq\bx$. Put
$n=\rk\bx$ and $\by=\pi^n_0(\bx)$. Let $i\in X$. {}From
$\bc_i\leq\bx$ and $\bc_i\in\cS(\Lambda)$ it follows that $\bc_i\leq\by$.
This holds for all $i\in X$, hence, as $\bc_{\xi}$ is clearly the
supremum of $\setm{\bc_i}{i\in X}$ in $\cS(\Lambda)$, we obtain that
$\bc_{\xi}\leq\by$. Therefore, $\bc_{\xi}\leq\bx$.
\end{proof}

Now we can state the main technical result of the paper. It says that
$\seqm{\bc_{\xi}}{\xi<\omega_1}$ is the least
non-eventually constant isotone $\omega_1$-sequence in $\cF(\omega_1)$
modulo the closed unbounded filter on $\omega_1$.

\begin{theorem}\label{T:NoIncrSeq}
Let $\sigma=\seqm{\bx_{\xi}}{\xi<\omega_1}$ be an isotone
$\omega_1$-sequence of elements of $\cF(\omega_1)$. Then either $\sigma$
is eventually constant or there exists a closed unbounded subset $C$ of
$\omega_1$ such that $\bc_{\xi}\leq\bx_{\xi}$ for all $\xi\in C$.
\end{theorem}

\begin{proof}
Assume that $\sigma$ is not eventually constant.
We put $X_{\xi}=\supp(\bx_{\xi})$ and
$n_{\xi}=|X_{\xi}|$, for all $\xi<\omega_1$. So
$\bx_{\xi}=\bx'_{\xi}[X_{\xi}/n_{\xi}]$, for some
$\bx'_{\xi}\in\cF(n_{\xi})$. As all sets~$X_{\xi}$ are finite, it
follows from the $\Delta$-Lemma (see
\cite[Lemma~22.6]{Jech}) that there are an uncountable subset~$I$ of
$\omega_1$ and a finite subset~$X$ of~$\omega_1$ such that
$X_{\xi}\cap X_{\eta}=X$ for all distinct $\xi,\eta\in I$. We may further
assume without loss of generality that there are $n<\omega$ and
$\bx\in\cF(n)$ such that $n_{\xi}=n$ and $\bx'_{\xi}=\bx$, for all
$\xi\in I$. Hence $\bx_{\xi}=\bx[X_{\xi}/n]$, for all $\xi\in I$. As
$\sigma$ is isotone but not eventually constant, it follows that $X$ is
a proper subset of $X_{\xi}$, for all $\xi\in I$. Put
$Y_{\xi}=X_{\xi}\setminus X$. Define $\rho(\xi)$ as the least element of
$Y_{\xi}$.

For subsets $U$ and $V$ of $\omega_1$, let $U<V$ hold, if $u<v$ for all
$\seq{u,v}\in U\times V$. By further shrinking $I$, we might assume that
$X<Y_{\xi}$, for all $\xi\in I$. In particular, observe
that $X=X_{\xi}\cap X_{\eta}$ is a lower subset of both $X_{\xi}$ and
$X_{\eta}$, for all $\xi\neq\eta$ in $I$.

Let $\xi<\eta$ in $I$ and suppose that there exists $\bz\in\cF(X)$ such
that $\bx_{\xi}\leq\bz\leq\nobreak\bx_{\eta}$. Applying the embedding
$\cF(e_{X_{\xi},X_{\eta}})$ to the inequality $\bx[X_{\xi}/n]\leq\bz$ and
using Lemma~\ref{L:Fix[Y/X]}, we obtain the inequality
$\bx[X_{\eta}/n]\leq\bz$, so $\bx_{\eta}=\bx[X_{\eta}/n]=\bz$, a
contradiction since the left hand side has support $X_{\eta}$ while the
right hand side has the smaller support $X$. Therefore, as
$\bx_{\xi}\leq\bx_{\eta}$ and by Lemma~\ref{L:Interp}, we obtain the
inequality $\bc_{\rho(\eta)}\leq\bx_{\eta}$.

Hence, we may assume that $\bc_{\rho(\xi)}\leq\bx_{\xi}$ for all
$\xi\in I$. It follows that
 \begin{equation}\label{Eq:Ineqrhobar}
 \bc_{\ol{\rho}(\xi)}\leq\bx_{\xi},\qquad\text{for all }\xi<\omega_1,
 \end{equation}
where we put
 \[
 \ol{\rho}(\xi)=\bigvee\famm{\rho(\eta)}{\eta\in I,\ \eta<\xi},
 \qquad\text{for all }\xi<\omega_1.
 \]
As the range of $\rho$ is unbounded, so is the range of $\ol{\rho}$.
Hence, as~$\ol{\rho}$ is a complete \jh\ from
$\omega_1$ to $\omega_1$, the set
$C=\setm{\xi<\omega_1}{\ol{\rho}(\xi)=\xi}$ is a closed unbounded subset
of~$\omega_1$. It follows from \eqref{Eq:Ineqrhobar} that the inequality $\bc_{\xi}\leq\bx_{\xi}$ holds for all $\xi\in C$.
\end{proof}

The following corollary expresses that $\cF(\omega_1)$ fails a certain
``monotone refinement property''.

\begin{corollary}\label{C:NoIncrSeq1}
There are no positive integer $n$ and no finite collection of isotone
$\omega_1$-sequences $\seqm{\bx_{i,\xi}}{\xi<\omega_1}$ of elements of
$\cF(\omega_1)$, for $0\leq i\leq n$, such that
\begin{enumerate}
\item $\bx_{0,\xi}=0$ and $\bx_{n,\xi}=\bc_{\xi}$, for all large enough
$\xi<\omega_1$.

\item $\bx_{i,\xi}\leq\bc_{\xi}$, for all $i\leq n$ and all large enough
$\xi<\omega_1$.

\item Either $\bx_{i+1,\xi}\leq\ba\vee\bx_{i,\xi}$ or
$\bx_{i+1,\xi}\leq\bb\vee\bx_{i,\xi}$, for all $i<n$ and all
$\xi<\omega_1$.
\end{enumerate}
\end{corollary}

\begin{proof}
We prove that for all $i\leq n$, there exists
$\eta_i<\omega_1$ such that $\bx_{i,\xi}\leq\bc_{\eta_i}$ for all
$\xi<\omega_1$. We argue by induction on~$i$. For $i=0$ it holds by
assumption, with $\eta_0=0$. Suppose that $\bx_{i,\xi}\leq\bc_{\eta_i}$,
for all $\xi<\omega_1$. Let $\xi>\eta_i$. Assume, for example, that
$\bx_{i+1,\xi}\leq\ba\vee\bx_{i,\xi}$; so
$\bx_{i+1,\xi}\leq\ba\vee\bc_{\eta_i}$. Observing that
$\bc_{\xi}\nleq\ba\vee\bc_{\eta_i}$, we get that
$\bc_{\xi}\nleq\bx_{i+1,\xi}$. As this holds for all
$\xi>\eta_i$ and by Theorem~\ref{T:NoIncrSeq}, we obtain that
$\seqm{\bx_{i+1,\xi}}{\xi<\omega_1}$ is eventually constant, and hence, by~(2), below some $\bc_{\eta_{i+1}}$, therefore completing the induction step.

In particular, for $i=n$, we obtain that the $\omega_1$-sequence
$\seqm{\bc_{\xi}}{\xi<\omega_1}$ is eventually dominated
by the constant~$\bc_{\eta_n}$, a contradiction.
\end{proof}

Hence we get a negative extension property for posets.

\begin{corollary}\label{C:NoIncrSeq2}
There are a poset measure
$\mu\colon(\omega_1+1)\times(\omega_1+1)\to\cF(\omega_1)$ such that
$\mu(\omega_1,0)=\ba\vee\bb$ but there are no poset $P$ containing
$\omega_1+1$, no poset measure
$\ol{\mu}\colon P\times P\to\cF(\omega_1)$ extending $\mu$, no positive
integer $n$, and no decomposition
$0=z_0\leq z_1\leq\cdots\leq z_n=\omega_1$ in~$P$ such that either
$\ol{\mu}(z_{i+1},z_i)\leq\ba$ or $\ol{\mu}(z_{i+1},z_i)\leq\bb$ for all
$i<n$.
\end{corollary}

\begin{proof}
Define $\mu\colon(\omega_1+1)\times(\omega_1+1)\to\cF(\omega_1)$ by
 \[
 \mu(\xi,\eta)=\begin{cases}
 0,&\text{if }\xi\leq\eta,\\
 \bc_{\xi},&\text{if }\eta<\xi<\omega_1,\\
 \ba\vee\bb,&\text{if }\eta<\xi=\omega_1.
 \end{cases}
 \qquad\text{for all }\xi,\eta\leq\omega_1.
 \]
It is straightforward to verify that $\mu$ is a poset measure on
$\omega_1+1$. Suppose that~$P$, $\ol{\mu}$, $n$, $z_0$, \dots, $z_n$
satisfy the given conditions. We put $\bx_{i,\xi}=\ol{\mu}(\xi,z_{n-i})$,
for all $\xi<\omega_1$. It is not hard, using the triangular inequality,
to verify that the elements $\bx_{i,\xi}$ satisfy the assumptions
(1)--(3) of Corollary~\ref{C:NoIncrSeq1}, a contradiction.
\end{proof}

As the following result shows, more ``amenable'' semilattices do satisfy a certain ``monotone refinement property''.

\begin{proposition}\label{P:MonRef}
Let $S$ be a distributive \jzs. If $S$ is either a lattice, or strongly
distributive, or countable, then for all $\ba,\bb\in S$, every chain
$\Lambda$, and every isotone $\Lambda$-sequence
$\seqm{\bc_i}{i\in\Lambda}$ of elements of $S$ such that
$\bc_i\leq\ba\vee\bb$ for all $i\in\Lambda$, then there are isotone
$\Lambda$-sequences $\seqm{\ba_i}{i\in\Lambda}$ and
$\seqm{\bb_i}{i\in\Lambda}$ of elements of $S$ such that $\ba_i\leq\ba$,
$\bb_i\leq\bb$, and $\bc_i=\ba_i\vee\bb_i$, for all $i\in\Lambda$.
\end{proposition}

\begin{proof}
If $S$ is a lattice the conclusion is trivial: put $\ba_i=\ba\wedge\bc_i$
and $\bb_i=\bb\wedge\bc_i$, for all $i\in\Lambda$.

Now assume that $S$ is strongly distributive. Denote by $C_i$ the
(finite) set of all maximal \jirr\ elements of $S$ below $\bc_i$, for all
$i\in\Lambda$. Observe that $C_i\subseteq\dnw C_j$, for all
$i\leq j$ in $\Lambda$.
For every finite subset $I$ of $\Lambda$, denote by $X_I$
the set of all families $\seqm{\seq{A_i,B_i}}{i\in I}$ such that
\begin{enumerate}
\item $A_i\subseteq\dnw\ba$, $B_i\subseteq\dnw\bb$, and $A_i\cup B_i=C_i$,
for all $i\in I$.

\item For all $i\leq j$ in $I$, $A_i\subseteq\dnw A_j$ and
$B_i\subseteq\dnw B_j$.
\end{enumerate}
We claim that $X_I$ is nonempty, for every finite subset $I$ of
$\Lambda$. We argue by induction on $|I|$. The conclusion is obvious for
$I=\es$. For $I=\set{i}$, put $A_i=\dnw\ba\cap C_i$ and
$B_i=\dnw\bb\cap C_i$. Now suppose that $I=\set{i}\cup J$, where $i<j$ for
all $j\in J$ and $J$ is nonempty. By induction hypothesis, there exists
an element $\seqm{\seq{A_k,B_k}}{k\in J}$ in~$X_J$. Put $j=\min J$,
$A_i=\dnw A_j\cap C_i$, and $B_i=\dnw B_j\cap C_i$. It is straightforward
to verify that $\seqm{\seq{A_k,B_k}}{i\in I}$ belongs to $X_I$. This
completes the induction step.

It follows that the set $\Omega_I$ of all families
$\seqm{\seq{A_i,B_i}}{i\in\Lambda}$ of elements of the Cartesian product
$\Omega=\prod\famm{\Pow(C_i)\times\Pow(C_i)}{i\in\Lambda}$
(where $\Pow(X)$ denotes the powerset of a set $X$) whose
restriction to $I$ belongs to $X_I$ is nonempty, for every finite
subset~$I$ of~$\Lambda$. Endow $\Omega$ with the product topology of the
discrete topologies on all (finite) sets $\Pow(C_i)\times\Pow(C_i)$. By
Tychonoff's Theorem, $\Omega$ is compact. Hence the intersection of all
$\Omega_I$, for $I$ a finite subset of $\Lambda$, is nonempty. Let
$\seqm{\seq{A_i,B_i}}{i\in\Lambda}$ be an element of that intersection.
Then the collection of all elements $\ba_i=\bigvee A_i$ and
$\bb_i=\bigvee B_i$, for $i\in\Lambda$, satisfies the required conditions.

Assume, finally, that $S$ is countable. Define an equivalence relation
$\equiv$ on $\Lambda$ by $i\equiv j$ if{f} $\bc_i=\bc_j$, for all
$i,j\in\Lambda$, and denote by $[i]$ the $\equiv$-equivalence class of
$i$, for any $i\in\Lambda$. Putting $\bc_{[i]}=\bc_i$ makes it possible
to replace $\Lambda$ by $\Lambda/{\equiv}$. In particular, as~$S$ is
countable, $\Lambda$ becomes countable as well. Now write
$\Lambda=\bigcup\famm{\Lambda_n}{n<\omega}$, where
$\seqm{\Lambda_n}{n<\omega}$ is an increasing sequence of finite subsets
of $\Lambda$ with $|\Lambda_n|=n$, for all $n<\omega$. Denote by $Y_n$
the set of all families $\seqm{\seq{\ba_l,\bb_l}}{l\in\Lambda_n}$ such
that $\ba_i\leq\ba$, $\bb_i\leq\bb$, $\bc_i=\ba_i\vee\bb_i$, and $i\leq j$
implies that $\ba_i\leq\ba_j$ and $\bb_i\leq\bb_j$, for all $i\leq j$ in
$\Lambda_n$. Suppose that we are given an element of $Y_n$ as above, and
denote by $k$ the unique element of $Y_{n+1}\setminus Y_n$. Suppose, for
example, that $\min\Lambda_n<k<\max\Lambda_n$, and denote by $i$ (resp.,
$j$) the largest (resp., least) element of $\Lambda_n$ below $k$ (resp.,
above $k$). As $\bc_k\leq\bc_j=\ba_j\vee\bb_j$, there are $\ba'\leq\ba_j$
and $\bb'\leq\bb_j$ such that $\bc_k=\ba'\vee\bb'$. Put
$\ba_k=\ba_i\vee\ba'$ and $\bb_k=\bb_i\vee\bb'$. Then
$\seqm{\seq{\ba_l,\bb_l}}{l\in\Lambda_{n+1}}$ is an element of $Y_{n+1}$.
So every element of~$Y_n$ extends to an element
of~$Y_{n+1}$. The proof is even easier in case either $k<\min\Lambda_n$ or
$k>\max\Lambda_n$. Hence we have constructed inductively a family
$\seqm{\seq{\ba_i,\bb_i}}{i\in\Lambda}$ whose restriction to $\Lambda_n$
belongs to $Y_n$, for all $n<\omega$. Therefore, the elements $\ba_i$ and
$\bb_i$, for $i\in\Lambda$, are as required.
\end{proof}

The ``monotone refinement property'' described above fails in $\cF(\omega_1)$. Indeed, consider the isotone $\omega_1$-sequence
$\seqm{\bc_{\xi}}{\xi<\omega_1}$ together with the inequalities
$\bc_{\xi}\leq\ba\vee\bb$, for $\xi<\omega_1$. Suppose that $\seqm{\ba_{\xi}}{\xi<\omega_1}$ and $\seqm{\bb_{\xi}}{\xi<\omega_1}$ are isotone $\omega_1$-sequences in $\cF(\omega_1)$ such that $\bc_{\xi}=\ba_{\xi}\vee\bb_{\xi}$ while $\ba_{\xi}\leq\ba$ and $\bb_{\xi}\leq\bb$, for all $\xi<\omega_1$.
Set $\bx_{0,\xi}=0$, $\bx_{1,\xi}=\ba_{\xi}$, and $\bx_{2,\xi}=\bc_{\xi}$, for all $\xi<\omega_1$. Then the isotone $\omega_1$-sequences $\seqm{\bx_{i,\xi}}{\xi<\omega_1}$, for $i\in\set{0,1,2}$, satisfy (1)--(3) of Corollary~\ref{C:NoIncrSeq1} (with $n=2$), a contradiction.

Nevertheless we do not know whether the monotone refinement property of~$S$ either implies or is implied by the statement that every $S$-valued poset measure extends to a V-measure.

\end{document}